\def\Ann{\mathop{\rm Ann}\nolimits}

\def\Ass{\mathop{\rm Ass}\nolimits}
\def\Assh{\mathop{\rm Assh}\nolimits}
\def\Att{\mathop{\rm Att}\nolimits}
\def\cd{\mathop{\rm cd}\nolimits}

\def\height{\mathop{\rm height}\nolimits}
\def\Hom{\mathop{\rm Hom}\nolimits}

\def\im{\mathop{\rm im}\nolimits}
\def\ker{\mathop{\rm ker}\nolimits}

\def\InjH{\mathop{\rm E}\nolimits}
\def\LCMo{\mathop{\rm H}\nolimits}

\def\Spec{\mathop{\rm Spec}\nolimits}

\def\Naturalsign{{\rm l\kern-.23em N}}
\parindent=0pt
\font \normal=cmr10 scaled \magstep0 \font \mittel=cmr10 scaled
\magstep1 \font \gross=cmr10 scaled \magstep5
\input amssym.def
\input amssym.tex
\gross Attached primes and Matlis duals of local cohomology modules
\normal
\bigskip
\bigskip
Michael Hellus \smallskip
 E-mail: michael.hellus@math.uni-leipzig.de
\smallskip
Universit\"at Leipzig, Fakult\"at f\"ur Mathematik und Informatik,
Mathematisches Institut, Ausgustusplatz 10/11, 04109 Leipzig,
Germany, Tel (+49) 341 9732186, Fax (+49) 341 9732199 \smallskip MSC
(2000): 13D45, 13E05, 13E10
\bigskip
\bigskip
{\bf Abstract} Let $J$ be an ideal of a noetherian local ring $R$.
We show new results on the set of attached primes $\Att _R(\LCMo
^l_J(R))$ of a local cohomology module $\LCMo ^l_J(R)$. To prove our
results we establish and use new relations between the set $\Att
_R(\LCMo ^l_J(R))$ of attached primes of a local cohomology module
and the set $\Ass _R(D(\LCMo ^l_J(R)))$ of associated primes of the
Matlis dual of the same local cohomology module.
\bigskip
The notions of attached primes and secondary decomposition of a
module were developed by MacDonald [15]. Attached primes of local
cohomology modules have been studied by MacDonald and Sharp. They
proved ([16, theorem 2.2])
$$\Att _R(\LCMo ^{\dim (M)}_\goth m(M))=\{ \goth p\in \Ass _R(M)\vert \dim
(R/\goth p)=\dim M\} $$ for every finite $R$-module $M$, where $R$
is a noetherian local ring. Dibaei and Yassemi generalized ([5,
theorem A]) this result to
$$\Att _R(\LCMo ^{\dim (M)}_\goth a(M))=\{ \goth p\in \Ass _R(M)\vert
\cd (\goth a,R/\goth p)=\dim (M)\} $$ for every finite $R$-module
$M$ and every ideal $\goth a$ of $R$. Besides reproving this result
(we remark that our new proof of this result is conceptual to some
degree), our main results are
$$\Att _R(\LCMo ^{n-1}_J(R))=\{ \goth p\in \Spec (R)\vert \dim
(R/\goth p)=n-1, \sqrt {\goth p+J}=\goth m\} \cup \Assh (R)$$ for
every $n$-dimensional local complete ring and any ideal $J$ of $R$
such that $\dim (R/J)=1$, $\LCMo ^n_J(R)=0$ (theorem 2.3 (i)) and
$$\{ \goth p\in \Spec (R)\vert x_1,\dots ,x_i\hbox { is a part of
system of parameters of }R/\goth p\} \subseteq \Att _R(\LCMo
^i_{(x_1,\dots ,x_i)R}(R))$$ for every local ring $(R,\goth m)$ and
any $x_1,\dots ,x_i\in R$ (theorem 2.3 (ii)). From these two main
results we draw some consequences on primes attached to local
cohomology modules ((2.4)--(2.6)).
\par
Our method to prove these results is to first establish some
elementary relations between attached and associated primes and then
to make use of what is known about associated primes of Matlis duals
of local cohomology (Hellus and St\"uckrad, [11], [9], [12]). This
method will also lead in a natural way to a proof of the
above-quoted theorem A from Dibaei and Yassemi.
\par
We need a careful (yet easy) analysis of attached primes, associated
primes and their relations; this program is carried out in section
1.\par Section 2 is devoted to results on the set of attached primes
of local cohomology modules, both known ones (theorem 2.1) and new
ones ((2.3) -- (2.5)). In addition we present new evidence for
conjecture (*) from [11, section 1]; this conjecture (*) says that
for every local ring $(R,\goth m)$ one has the equality
$$\Ass _R(D(\LCMo ^i_{(x_1,\dots ,x_i)R}(R)))=\{ \goth p\in \Spec
(R)\vert \LCMo ^i_{(x_1,\dots ,x_i)R}(R/\goth p)\neq 0\} $$ for any
$x_1,\dots ,x_i\in R$, where $D$ stands for the Matlis duality
functor with respect to any fixed $R$-injective hull of $R/\goth m$
(see [3], [4], [6], [17] for details on Matlis duality and injective
modules and [11, section 1] for comments on conjecture (*)). The
details of this new evidence are in (2.5) and (2.6), the basic idea
is simple: First, we know from [11] that conjecture (*) is
equivalent to the fact that $\Ass _R(D(\LCMo ^i_{(x_1,\dots
,x_i)R}(R)))=:X$ is closed under generalization, i. e. $\Spec (R)\ni
\goth p_0 \subseteq \goth p_1\in X$ implies $\goth p_0\in X$.
Second, in section 1 some relations between $X$ and $\Att _R(\LCMo
^i_{(x_1,\dots ,x_i)R}(R))$ are established. And third, theorem
(1.12) shows that a property which is a consequence of stableness
under generalization holds for $\Att _R(\LCMo ^i_{(x_1,\dots
,x_i)R}(R))$. The combination of these three facts will lead to new
evidence for conjecture (*) in a special case.
\bigskip \bigskip \mittel 1.
Notions \normal
\bigskip
Associated primes are clearly related to the notion of primary
decomposition; in a dual way, attached primes are related to
secondary decomposition. The theories of primary resp. secondary
decomposition are well-known, information on it can be found in many
textbooks (e. g. [6], [18]). But as we need quite general statements
where the module in question is not necessarily finite or artinian,
we carefully repeat what is known; in this section we omit proofs as
they are the same like in the finite resp. artinian case.
\bigskip (1.1)
Let $R$ be a ring, $M\neq 0$ an $R$-module and $N$ an $R$-submodule
of $M$. We say $M$ is coprimary iff the following condition holds:
For every $x\in R$ the endomorphism $M\buildrel x\over \to M$ given
by multiplication by $x$ is injective or nilpotent (i. e. $\exists
N\in \Naturalsign : x^N\cdot M=0$). If $M$ is coprimary $\sqrt {
\Ann _R(M)}$ is a prime ideal of $R$. In general we say $N$ is a
primary submodule of $M$ iff $M/N$ is coprimary. Now let $U_1,\dots
,U_s\subseteq M$ be submodules of $M$. We say the $s$-tuple $(U_1,
\dots ,U_s)$ is a primary decomposition of (the zero ideal of) $M$
iff the following two conditions hold:
\par
(i) $U_1\cap \dots \cap U_s=0$.
\par
(ii) All $U_i$ are primary submodules of $M$.
\smallskip
In this case $(U_1, \dots ,U_s)$ is called a minimal primary
decomposition of $M$ iff, in addition, the following two statements
hold:
\par
(iii) Every $U_1\cap \dots \cap \hat {U_i}\cap \dots \cap U_s$ is
not zero.
\par
(iv) The ideals $\sqrt {\Ann _R(M/U_i)}$ (for $i=1,\dots ,s$) are
pairwise different.
\smallskip
It is clear that if there exists a primary decomposition of $M$
there is also a minimal one.
\smallskip
(1.2) Let $R$ be a noetherian ring, $M$ an $R$-module and assume
there exists a minimal primary decomposition $(U_1,\dots ,U_s)$ of
$M$. Then the set
$$\{ \sqrt {\Ann _R(M/U_i)}\vert i=1,\dots ,s\} =: \Ass _R(M)$$
does not depend on the choice of a minimal primary decomposition of
$M$ (the proof of this goes just like the well-known proof in case
$M$ is finite). We say the prime ideals of $\Ass _R(M)$ are
associated to $M$
\smallskip
(1.3) Let $R$ be a noetherian ring and $M$ a noetherian (i. e.
finitely generated) $R$-module. Then it is well-known that $M$ has a
(minimal) primary decomposition. Note that this holds without the
hypothesis $R$ is noetherian, but anyway $M$ being noetherian
implies that $R/\Ann _R(M)=:\overline R$ is noetherian and $M$ is a
$\overline R$-module.
\smallskip
(1.4) Let $R$ be a noetherian ring and $M$ an $R$-module. One
defines
$$\Ass _R(M):=\{ \goth p\subseteq R\hbox { prime ideal }\vert
\exists m\in M:\goth p=\Ann _R(m)\} .$$ It is easy to see that this
definition agrees with the above one whenever $M$ has a primary
decomposition.
\smallskip
(1.5) Let $R$ be a ring and $M\neq 0$ an $R$-module. One says $M$ is
secondary iff for every $x\in R$ the endomorphism $M\buildrel x\over
\to M$ given by multiplication by $x$ is either surjective or
nilpotent. Now let $M$ be arbitrary and $U_1,\dots, U_s\subseteq M$
$R$-submodules. We say the $s$-tuple $(U_1,\dots ,U_s)$ is a
secondary decomposition of $M$ iff the following two conditions
hold: $U_1+\dots +U_s=M$ and all $U_i$ are secondary. In this case
the secondary decomposition $(U_1,\dots ,U_s)$ is called minimal iff
the following two conditions hold: All $U_1+\dots +\hat {U_i}+\dots
+U_s$ are proper subsets of $M$ and all $\sqrt {\Ann _R(U_i)}$ are
pairwise different. Again existence of a secondary decomposition
implies existence of a minimal one.
\smallskip
(1.6) Let $R$ be a noetherian ring and $M$ an $R$-module; assume
there exists a minimal secondary decomposition $(U_1,\dots ,U_s)$ of
$M$. Then the set
$$\Att _R(M):=\{ \Ann _R(U_i)\vert i=1,\dots s\} $$
does not depend on the choice of a minimal secondary decomposition
of $M$. We say the prime ideals in $\Att _R(M)$ are attached to $M$.
\smallskip
(1.7) Let $R$ be a noetherian ring and $M$ an artinian $R$-module.
Then there exists a (minimal) secondary decomposition of $M$. The
proof is simply a dual version of the proof of (1.3). Again this
works also if $R$ is not noetherian.
\smallskip
(1.8) Let $R$ be a noetherian ring and $M$ an $R$-module. We define
$$\Att _R(M):= \{ \goth p\subseteq R \hbox { prime ideal }\vert
\exists \hbox { an }R\hbox {-submodule }U \subseteq M: \goth p=\Ann
_R(M/U)\} .$$ Is is not very difficult to see that this definition
agrees with the first one if $M$ has a secondary decomposition.
\smallskip
(1.9) Let $(R,\goth m)$ be a noetherian local ring, $M$ an
$R$-module and $(U_1,\dots ,U_s)$ a minimal primary decomposition of
$M$. By $D(\_ )$ we denote the Matlis dual functor from the category
of $R$-modules to itself sending $M$ to $\Hom _R(M,\InjH )$, where
$\InjH :=\InjH _R(R/\goth m)$ shall denote an $R$-injective hull of
$R/\goth m$. The following implications are clear by duality:
\smallskip
(i) $U_1\cap \dots \cap U_s=0 \Rightarrow D(M/U_1)+\dots +\dots
D(M/U_s)=D(M)$
\par
(ii) $M/U_i$ is coprimary $\Rightarrow D(M/U_i)$ is secondary (for
every $i$)
\par
(iii) The primary decomposition $(U_1,\dots ,U_s)$ of $M$ is minimal
$\Rightarrow $\par
the secondary decomposition $(D(M/U_1),\dots
,D(M/U_s))$ of $D(M)$ is minimal.
\par
(iv) $\Ann _R(M/U_i)=\Ann _R(D(M/U_i))$ (for every $i$)
\smallskip
Thus we have
$$\Ass _R(M)=\Att _R(D(M))\ \ .$$
In a very similar way the
following statement holds: Any (minimal) secondary decomposition of
$M$ induces a (minimal) primary decomposition of $D(M)$. In
particular, if $M$ has a secondary decomposition:
$$ \Att _R(M)=\Ass _R(D(M))\ \ .$$
Remark: It is true that if $U_1,\dots ,U_s$ are arbitrary submodules
of $R$ such that $(D(M/U_1),\dots ,D(M/Us))$ is a (minimal)
secondary decomposition of $D(M)$ then $(U_1,\dots ,U_s)$ is a
(minimal) primary decomposition of $M$, but note that we do not know
that every submodule of $D(M)$ is of the form $D(M/U)$ for some
submodule $U$ of $M$. Similarly, if $U_1,\dots ,U_s$ are arbitrary
submodules of $M$ such that $(D(M/U_1),\dots ,D(M/U_s))$ is a
(minimal) primary decomposition of $D(M)$ then $(U_1,\dots ,U_s)$ is
a (minimal) secondary decomposition of $M$.
\smallskip
(1.10) Let $(R,\goth m)$ be a noetherian local ring, $\goth p$ a
prime ideal of $R$ and $M$ an $R$-module. Then
$$\eqalign {\goth p\in \Ass _R(M)&\iff \exists \hbox { finitely generated
submodule }U \hbox { of } M: \goth p=\Ann _R(U),\cr \goth p\in \Att
_R(D(M))&\iff \exists \hbox { submodule }U^\prime \hbox { of }D(M):
\goth p=\Ann _R(D(M)/U^\prime ).\cr }$$ In particular the existence
of a submodule $U$ of $M$ satisfying $\goth p=\Ann _R(U)$ implies
$\goth p\in \Att _R(D(M))$. Therefore we have
$$\Ass _R(M)\subseteq \Att _R(D(M)).$$
This inclusion is strict in general: Take for example $M=\InjH
=\InjH _R(R/\goth m)$, an $R$-injective hull of $R/\goth m$: $\Ass
_R(\InjH _R(R/\goth m))=\{ \goth m\} $, but $D(\InjH )=\hat R$ and
so $\Att _R(D(\InjH ))=\Spec (R)$. But nevertheless a stronger
inclusion holds (plug in $D(M)$ for $M$ to see that it is actually
stronger):
\bigskip
(1.11) {\bf Theorem}
\par
Let $(R,\goth m)$ be a noetherian local ring and $M$ an $R$-module.
Then
$$\Ass _R(D(M))\subseteq \Att _R(M)$$
and the sets of prime ideals maximal in each side respectively
coincide:
\par
$$\{ \goth p\vert \goth p\hbox { maximal in } \Ass _R(D(M))\}=\{ \goth p\vert \goth p\hbox { maximal in } \Att
_R(M)\}.$$ Proof: \par Let $\goth p\in \Ass _R(D(M))$ be arbitrary.
There exists a submodule $U^\prime $ of $D(M)$ such that $U^\prime
=R\cdot u^\prime \cong R/\goth p$ for some $u^\prime \in U^\prime
\subseteq D(M)$. $u^\prime $ induces a monomorphism $\overline
{u^\prime }:M/\ker (u^\prime )\to \InjH $ and so we have
$$\goth p=\Ann _R(U^\prime )=\Ann _R(u^\prime )=\Ann _R(\overline
{u^\prime })=\Ann _R(M/\ker (u^\prime ));$$ this implies $\goth p\in
\Att _R(M)$. Having proved this we only have to show that an
arbitrary prime ideal $\goth p$ of $R$ which is maximal in $\Att
_R(M)$ is associated to $D(M)$: $\goth p\in \Att _R(M)$ implies
$M/\goth pM\neq 0$ and so we must have $\Hom _R(R/\goth
p,D(M))=D(M/\goth pM)\neq 0$; but by the maximality hypothesis on
$\goth p$ implies $\goth p\in \Ass _R(D(M))$.
\bigskip
(1.12) {\bf Theorem}
\par
Let $(R,\goth m)$ be a noetherian local ring and $M$ an $R$-module.
Assume $(\goth p_i)_{i\in \Naturalsign }$ is a sequence of prime
ideals attached to $M$; assume furthermore that $\goth q:=\bigcap
_{i\in \Naturalsign }\goth p_i$ is a prime ideal of $R$. Then $\goth
q$ is also attached to $M$.
\par
Proof: \par For every $i$ we choose a quotient $M_i$ of $M$ such
that $\Ann _R(M_i)=\goth q_i$. Now the canonically induced map
$\iota :M\to \prod _{i\in \Naturalsign }M_i$ induces a surjection
$M\to \im (\iota )$; we obviously have $\bigcap _{i\in \Naturalsign
}\goth p_i\subseteq \Ann _R(\im (\iota ))$; on the other hand, for
every $i$ and every $s\in R\setminus \goth p_i$ there is a
$\overline {m_i}\in M_i$ coming from an element $m_i\in M$ that has
$s\cdot \overline {m_i}\neq 0$. But this implies that $s$ cannot
annihilate $\im (\iota )$; therefore $\Ann _R(\im (\iota))=\bigcap
_{i\in \Naturalsign }\goth p_i=\goth q$ and the statement follows.
\bigskip
\bigskip
\mittel 2. Results \normal
\bigskip
There are some results on the set of attached primes of local
cohomology modules: In [16, theorem 2.2] it was shown that if
$(R,\goth m)$ is a noetherian local ring and $M$ is a finitely
generated $R$-module then
$$\Att _R(\LCMo ^{\dim (M)}_\goth m(M))=\{ \goth p\in \Ass _R(M)\vert
\dim (R/\goth p)=\dim (M)\} .$$ In [5, Theorem A] this was
generalized to
$$\Att _R(\LCMo ^{\dim (M)}_\goth a(M))=\{ \goth p\in \Ass
_R(M)\vert \cd (\goth a,R/\goth p)=\dim (M)\} ,$$ where $\goth
a\subseteq R$ is an ideal and $\cd (\goth a,R/\goth p):=\max \{ l\in
\Naturalsign \vert \LCMo ^l_\goth a(R/\goth p)\neq 0\} $. We are
going to show (theorem 2.1) that the results of section 1 lead to a
natural proof of this theorem and, furthermore, to new results on
the attached primes of local cohomology modules ((2.3) -- (2.6)).
\bigskip
Let $(R,\goth m)$ be a noetherian local
$n$-dimensional ring and $\goth a\subseteq R$ an ideal. Then $\LCMo
^n_\goth a(R)$ is an artinian $R$-module and hence
$$\Ass _R(D(\LCMo ^n_\goth a(R)))=\Att _R(\LCMo ^n_\goth a(R)).$$
Now assume that we have ($\LCMo ^n_\goth a(R)\neq 0$ and) $\goth
p\in \Att _R(\LCMo ^n_\goth a(R))$; then we get
$$0\neq \LCMo ^n_\goth a(R)/\goth p\LCMo ^n_\goth a(R)=\LCMo
^n_\goth a(R/\goth p),$$ i. e. $\goth p\in \Assh (R) $($:=\{ \goth
q\in \Spec(R)\vert \dim (R/\goth q)=\dim (R)\}$) and $\cd (\goth
a,R/\goth p)=n$.
\par
Now suppose conversely that we have a prime ideal $\goth p$ of $R$
such that $\cd (\goth a,R/\goth p)=n$, equivalently $\LCMo ^n_\goth
a(\hat {R/\goth p})\neq 0$. By Hartshorne-Lichtenbaum vanishing we
get a prime ideal $\goth q\subseteq \hat R$ satisfying $\goth
p=\goth q\cap R$ and $\sqrt {\goth a\hat R+\goth q}=\goth m_{\hat
R}$($:=$maximal ideal of $\hat R$); this in turn implies
$$0\neq \LCMo ^n_{\goth a\hat R}(\hat R/\goth q)=\LCMo ^n_{\goth
m_{\hat R}}(\hat R/\goth q).$$ Matlis duality theory shows that
$\goth q\in \Ass _{\hat R}(D(\LCMo ^n_{\goth a\hat R}(\hat R)))$. It
is easy to see that
$$D(\LCMo ^n_{\goth a\hat R}(\hat R))=D(\LCMo ^n_\goth a(R)),$$
holds canonically, the $D$-functors taken over $\hat R$ resp. over
$R$. Thus we have shown
$$\Att _R(\LCMo ^n_\goth a(R))=\{ \goth p\vert \cd (\goth a,R/\goth
p)=n\} .$$ For every finitely generated $R$-module $M$ we can apply
this result to the ring $R/\Ann _R(M)$ and we get
\bigskip
(2.1) {\bf Theorem}
\par
Let $(R,\goth m)$ be a noetherian local ring and $M$ a finitely
generated $n$-dimensional $R$-module. Then
$$\Att _R(\LCMo ^n_\goth a(M))=\{ \goth p\in \Ass _R(M)\vert \cd
(\goth a,R/\goth p)=n\} .$$ (2.2) {\bf Remark} This is [5, Theorem
A], where it was proved by different means.
\bigskip
In section 1 we established several relations between attached
primes of a module and associated primes of the Matlis dual of the
same module; theorem 2.1 is a consequence of these relations; we can
retrieve more information out of these to get new theorems on the
attached primes of top local cohomology modules:
\bigskip
(2.3) {\bf Theorem} \par Let $(R,\goth m)$ be a $d$-dimensional
noetherian local ring.
\par
(i) If $J$ is an ideal of $R$ such that $\dim (R/J)=1$ and $\LCMo
^d_J(R)=0$ then
$$\Assh (R)\subseteq \Att _R(\LCMo ^{d-1}_J(R))\ \ .$$
If, in addition, $R$ is complete, one has
$$\Att _R(\LCMo ^{d-1}_J(R))=\{ \goth p\in \Spec (R)\vert \dim
(R/\goth p)=d-1, \sqrt {\goth p+J}=\goth m\} \cup \Assh (R).$$ \par
(ii) For any $x_1,\dots ,x_i\in R$ there is an inclusion
$$\{ \goth p\in \Spec (R)\vert x_1,\dots ,x_i \hbox { is a part of a
system of parameters of }R/\goth p\} \subseteq \Att _R(\LCMo
^i_{(x_1,\dots ,x_i)R}(R)).$$ Proof: \par (i) Note that [12,
theorems 5.4 and 5.5] show that one has $\Assh (R)=\Assh (D(\LCMo
^{d-1}_J(R)))$ in the given situation and, if $R$ is complete, $\Ass
_R(D(\LCMo ^{d-1}_J(R)))=\{ \goth p\in \Spec (R)\vert \dim (R/\goth
p)=1, \dim (R/(\goth p+J))=0\} \cup \Assh (R)$. Now we use theorem
(1.11)] and remark: If $R$ is complete, given an arbitrary $\goth
p\in \Att _R(\LCMo ^{d-1}_J(R))$ it follows that $\LCMo
^{d-1}_J(R/\goth p)\neq 0$ and hence, by Hartshorne Lichtenbaum
vanishing, that $\dim (R/\goth p)\geq d-1$ and, if $\dim (R/\goth
p)=d-1$, that $\goth p+J$ is $\goth m$-primary.
\par
(ii) Follows from theorem (1.11) and [12, theorem 1.3 (ii)].
\bigskip
(2.4) {\bf Corollary}
\par
For every $x\in R$ one has
$$\Att _R(\LCMo ^1_{xR}(R))=\Spec (R)\setminus {\frak V}(x).$$
Proof: \par "$\subseteq $" Let $\goth p\in \Att _R(\LCMo
^1_{xR}(R))$. Then $0\neq \LCMo ^1_{xR}(R)/\goth p\LCMo
^1_{xR}(R)=\LCMo ^1_{xR}(R/\goth p) \Rightarrow x\not\in \goth p$.
"$\supseteq $" follows from [12, theorem 1.3 (ii)].
\bigskip
(2.5) {\bf Remarks}
\par
(i) It was shown in [11, theorem 2.2.1] that for any
$x_1,\dots,x_i\in R$ there is an inclusion
$$\Ass _R(D(\LCMo ^i_{(x_1,\dots ,x_i)R}(R)))\subseteq \{ \goth p\in
\Spec (R)\vert \LCMo ^i_{(x_1,\dots ,x_i)R}(R/\goth p)\neq 0\}.$$ By
what we have proved so far it is clear that there is chain of
inclusions
$$\Ass _R(D(\LCMo ^i_{(x_1,\dots ,x_i)R}(R)))\subseteq \Att _R(\LCMo
^i_{(x_1,\dots ,x_i)R}(R))\subseteq \{ \goth p\in \Spec (R)\vert
\LCMo ^i_{(x_1,\dots ,x_i)R}(R/\goth p)\neq 0\} .$$ (ii) In [11,
section 1] it was conjectured that the inclusion
$$\Ass _R(D(\LCMo ^i_{(x_1,\dots ,x_i)R}(R)))\subseteq \{ \goth p\in
\Spec (R)\vert \LCMo ^i_{(x_1,\dots ,x_i)R}(R/\goth p)\neq 0\}$$ is
always an equality; we denote this conjecture by (*); if true, it
implies immediately
$$\Ass _R(D(\LCMo ^i_{(x_1,\dots ,x_i)R}))=\Att _R(\LCMo
^i_{(x_1,\dots ,x_i)R}(R)).$$ (iii) In the situation of theorem
(2.3) (i) the attached primes of the top local cohomology module
coincide with the associated primes of the Matlis dual of the top
local cohomology module.
\smallskip
(2.6) We now assume that $k$ is a field and $R=k[[X_1,\dots ,X_n]]$
is a power series algebra in $n$ variables $X_1,\dots ,X_n$; let
$i\in \{ 1,\dots ,n\} $. [12, section 3] and theorem 1.11 imply the
following statements:
\smallskip
(i) $i=n$: $\Att _R(\LCMo ^n_{(X_1,\dots ,X_n)R}(R)=\{ 0\} $.
\par
(ii) $i=n-1$: $\Att _R(\LCMo ^{n-1}_{(X_1,\dots ,X_{n-1})}(R)=\{ 0\}
\cup \{ pR\vert p\in R$ prime element, $p\not\in (X_1,\dots
,X_{n-1})R\} $.
\par
(iii) $i=n-2$: \par - $\{ 0\} \in \Att _R(\LCMo ^{n-2}_{(X_1,\dots
,X_{n-2})R}(R))$; \par - if $\goth p$ is a height 2 prime ideal of
$R$ such that $\sqrt {(X_1,\dots ,X_{n-2})R+\goth p}=\goth m$ then
$\goth p\in \Att _R(\LCMo ^{n-2}_{(X_1,\dots ,X_{n-2})R}(R))$;
\par
- conversely, $\goth p\in \Att _R(\LCMo ^{n-2}_{(X_1,\dots
,X_{n-2})R}(R))$ implies that $\height (\goth p)\leq 2$; \par - if
$p\in R$ is a prime element such that $p\not\in (X_1,\dots
,X_{n-2})R$ then $pR\in \Att _R(\LCMo ^{n-2}_{(X_1,\dots
,X_{n-2})R}(R))$;\par -
if $p\in R$ is a minimal generator of
$(X_1,\dots ,X_{n-2})R$ then $pR$ is not attached to $\LCMo
^{n-2}_{(X_1,\dots ,X_{n-2})R}(R)$.
\smallskip
In [12, 3.3 (i)] it was shown that for a prime element $p\in
(X_1,\dots ,X_{n-2})R$ under certain conditions there exist
infinitely many (pairwise different) prime ideals $(\goth p_l)_{l\in
\Naturalsign }$ of height 2 attached to $\LCMo ^{n-2}_{(X_1,\dots
,X_{n-2})R}(R)$ and containing $p$. As any $q\in \bigcap _{l\in
\Naturalsign }\goth p_l$ must satisfy $\height (p,q)R<2$ it is clear
that we have $pR=\bigcap _{l\in \Naturalsign }\goth p_l$. Now
theorem (1.12) implies $pR\in \Att _R(\LCMo ^{n-2}_{(X_1,\dots
,X_{n-2})R}(R))$. But in view of [12, theorem 1.1] and theorem
(1.11) it is clear that $pR\in \Att _R(\LCMo ^{n-2}_{(X_1,\dots
,X_{n-2})R}(R)$ is a necessary condition for conjecture (*). This
gives new evidence for conjecture (*).
 \vfil \eject {\bf References} \normal
\smallskip
\parindent=0.8cm
\def\litem{\par\noindent \hangindent=\parindent\ltextindent}
\def\ltextindent#1{\hbox to \hangindent{#1\hss}\ignorespaces}
\litem{1.} Bass, H. On the ubiquity of Gorenstein rings, {\it Math. Z.} {\bf
82}, (1963) 8--28.
\medskip
\litem{2.} Brodmann, M. and Hellus, M. Cohomological patterns of coherent
sheaves over projective schemes, {\it Journal of Pure and Applied Algebra}
{\bf 172}, (2002) 165--182.
\medskip
\litem{3.} Brodmann, M. P. and Sharp, R. J. Local Cohomology, {\it Cambridge
studies in advanced mathematics} {\bf 60}, (1998).
\medskip
\litem{4.} Bruns, W. and Herzog, J. Cohen-Macaulay Rings, {\it
Cambridge University Press}, (1993).
\medskip
\litem{5.} Dibaei, M. T. and Yassemi, S. Attached primes of the top
local cohomology modules with respect to an ideal, {\it Arch. Math.}
{\bf 84}, (2005) 292--297.
\medskip
\litem{6.} Eisenbud, D. Commutative Algebra with A View Toward
Algebraic Geometry, {\it Springer Verlag}, (1995).
\medskip
\litem{7.} Grothendieck, A. Local Cohomology, {\it Lecture Notes in
Mathematics, Springer Verlag}, (1967).
\medskip
\litem{8.} Hellus, M. Local Homology, Cohen-Macaulayness and
Cohen-Macaulayfications, to appear in Algebra Colloquium.
\medskip
\litem{9.} Hellus, M. Matlis duals of top local cohomology modules
and the arithmetic rank of an ideal, preprint.
\medskip
\litem{10.} Hellus, M. On the set of associated primes of a local
cohomology module, {\it J. Algebra} {\bf 237}, (2001) 406--419.
\medskip
\litem{11.} Hellus, M. On the associated primes of Matlis duals of
top local cohomology modules, to appear in {\it Communications in
Algebra} {\bf 33}.
\medskip
\litem{12.} Hellus, M. and St\"uckrad, J. Matlis duals of top Local
Cohomology Modules, preprint.
\medskip
\litem{13.} Huneke, C. Problems on Local Cohomology, {\it Res. Notes
Math. } {\bf 2}, (1992) 93--108.
\medskip
\litem{14.} Huneke, C. and Lyubeznik, G. On the vanishing of local
cohomology modules, {\it Invent. math.} {\bf 102}, (1990) 73--93.
\medskip
\litem{15.} MacDonald I. G. Secondary representation of modules over
a commutative ring, {\it Symp. Math.} {\bf XI}, (1973) 23--43.
\medskip
\litem{16.} MacDonald I. G. and Sharp, R. Y. An elementary proof of
the non-vanishing of certain local cohomology modules, {\it Quart.
J. Math. Oxford} {\bf 23}, (1972) 197--204.
\medskip
\litem{17.} Matlis, E. Injective modules over Noetherian rings, {\it
Pacific J. Math.} {\bf 8}, (1958) 511--528.
\medskip
\litem{18.} Matsumura, H. Commutative ring theory, {\it Cambridge
University Press}, (1986).
\medskip
\litem{19.} Scheja, G. and Storch, U. Regular Sequences and
Resultants, {\it AK Peters}, (2001).
\end